\documentclass[12pt,bezier]{article}

\usepackage{amsmath, amssymb, amsfonts, amsthm,url}
\usepackage{xcolor,colortbl}
\usepackage{graphics}
\usepackage{amssymb,amsmath,amsthm,latexsym,tikz,url,float,subfig,caption,graphicx,pgfplots,mathrsfs,hyperref}
 \usetikzlibrary{matrix}
\tikzstyle{vertex}=[circle, draw, inner sep=0pt, minimum size=3pt]
\newcommand{\vertex}{\node[vertex]}
\usepackage{color}
\newtheorem{theorem}{Theorem}
\newtheorem{lemma}[theorem]{Lemma}
\theoremstyle{definition}

\newtheorem{example}[theorem]{Example}

\newtheorem{remark}[theorem]{Remark}

\begin{document}

\title{Some results on perfect codes in Cayley sum graphs}

\author{ Masoumeh Koohestani$^{\,\rm a}$\quad Doost Ali Mojdeh  $^{\,\rm a}$ \quad\\
 Mohsen Ghasemi $^{\,\rm b}$ Hassan Khodaiemehr $^{\,\rm c}$
\\[.3cm]
{\sl\normalsize $^{\rm a}$ Department of Mathematics,  University of Mazandaran,}\\
{\sl\normalsize  P. O. Box 47416-1467, Babolsar, Iran}\\
{\sl\normalsize $^{\rm b}$Department of Mathematics, Urmia University,}\\
{\sl\normalsize  P. O. Box 575615-1818, Urmia, Iran}\\
{\sl\normalsize $^{\rm c}$ School of Engineering, Faculty of Applied Science, University of British Columbia (UBC), }\\
{\sl\normalsize  Kelowna, BC, Canada}
}
\maketitle
\footnote{{\em E-mail Addresses}: { m.koohestani@umail.umz.ac.ir, dmojdeh@umz.ac.ir, 
m.ghasemi@urmia.ac.ir, hassan.khodaiemehr@ubc.ca }}

\maketitle
\begin{abstract} 
This paper explores perfect codes within the context of Cayley sum graphs. We establish a necessary and sufficient condition for a normal subgroup of a finite group to be a subgroup perfect code and another necessary and sufficient condition for a normal subgroup not to be a subgroup perfect code. These results enable us to identify relationships between perfect codes of a group and its quotient groups. Furthermore, we examine a special nonabelian group, $V_{8n}$, and identify its possible connected Cayley sum graphs, as well as their corresponding subgroup perfect codes and subgroup total perfect codes.


 \vspace{4mm}
 	
\noindent {\bf AMS Mathematics Subject Classification\,(2020):}   05C25, 05C69, 94B99
\end{abstract}

\section{Introduction}
Perfect codes have played a pivotal role in information theory ever since the advent of coding theory in the late 1940s. Refer to surveys \cite{H, van} for a comprehensive collection of results on perfect codes. The significance of perfect codes is widely acknowledged, with Hamming and Golay codes being noteworthy examples. The concept of perfect codes can be extended to graphs in a natural manner \cite{B}. Within the context of graphs, perfect codes are also known as efficient dominating sets \cite{DS, KP} or independent perfect dominating sets \cite{L}.

Perfect codes in Cayley graphs have
 garnered significant attention, as evidenced by recent research in the field \cite{CWX, MWWZ, Zh, Zho}. For a group $G$ with the identity element 1 and an inverse-closed subset of $G \setminus \{1\}$ like $X$ (that is, $X^{ -1} := \{x^{ -1} ~|~ x \in X\} = X$), the \textit{Cayley graph} ${\rm{Cay(G, X)}}$ of $G$ with the \textit{connection set} $X$ is a graph whose vertex set is $G$ and two elements $x, y \in  G$ are adjacent if and only if $yx^{-1} \in X$. In \cite{HXZ}, Huang et al. introduced the following concept: A subset $C$ of a group $G$ is called
a \textit{perfect code} of $G$ if there exists a Cayley graph ${\rm Cay}(G, X)$ of $G$ which admits $C$ as a
perfect code.
 
 Of particular interest are perfect codes in Cayley sum graphs, which have drawn considerable attention from researchers. Cayley sum graphs, also referred to as addition graphs \cite{CGW}, addition Cayley graphs \cite{GLS, L, SGS}, and sum graphs \cite{C}, are variations of Cayley graphs. The concept of Cayley sum graphs, initially developed for abelian groups, was later generalized for arbitrary groups in \cite{AT}. Numerous authors have conducted research on Cayley sum graphs, contributing to the body of knowledge in this area (see \cite{MFW, MWY, Zhang-2024}). This paper builds upon this foundation by exploring perfect codes in Cayley sum graphs. 
 
Let $\Gamma$  be a graph with vertex set $V (\Gamma)$ and edge set $E(\Gamma)$.  A \textit{perfect code}
 in $\Gamma$ is a subset $C$ of $V(\Gamma)$ such that no two vertices in $C$ are adjacent and every vertex in $V(\Gamma)\setminus C$ is adjacent to exactly one vertex in $C$. Also, $C$ is called a \textit{total perfect code} in  $\Gamma$,  if every vertex in  $V (\Gamma)$ has exactly one neighbor in $C$. 
Let $G$ be a group. An element $x$ of $G$ is called a \textit{square} if $x = y^2$ for some element
$y \in G$. A subset of $G$ is called a  \textit{square-free subset} of $G$ if it is a set without squares. A
subset $X$ of $G$ is called a  \textit{normal subset} if $X$ is a union of some conjugacy classes of $G$ or
equivalently, for every $g \in G$, $g^{-1}Xg := \{g^{-1}xg : x \in X\} = X$. Remark that any subset of an
abelian group is normal. Let $X$ be a a normal subset of $G$. The  \textit{Cayley sum graph} ${\rm CS}(G, X)$ of $G$ with respect to the
 \textit{connection set} $X$ is a graph whose vertex set is $G$ and  two vertices $g$ and $h$ being
adjacent if and only if $gh \in X$ and $g \neq h$. 
A perfect code (total perfect code) of  ${\rm CS}(G, X)$ which is also a subgroup of $G$ is called
a  \textit{subgroup perfect code} (\textit{subgroup total perfect code}) of $G$. A subgroup perfect code (subgroup total perfect code) $H$  of the group $G$ is said to be \textit{nontrivial} when $H \neq \{1\}$. Also, the set of all perfect codes (total perfect codes) of ${\rm CS}(G, X)$ is denoted by $\mathcal{P} (G, X)$ ($\mathcal{T} (G, X)$). 
Additionally, it is obvious that the neighborhood
of a vertex $g$ is $Xg^{-1}$ if $g^2 \notin X$ and $(X \setminus \{g^2\})g^{-1}$ if $g^2 \in X$. 
Therefore ${\rm CS}(G, X)$ is
a regular graph if and only if either $X$ is square-free in $G$ or $X$ consists of all squares
of $G$. 

In what follows, we summarize the main  points that we will cover in this paper. In Section \ref{Results:SPC}, we present necessary and sufficient conditions for a subgroup of a finite group to be a subgroup perfect code, drawing interesting connections between subgroup perfect codes of a group and those of its quotient group. 
Section \ref{special} is dedicated to identifying perfect codes and total perfect codes within the renowned group $V_{8n}$.
 Detailed computational methods employed to define Cayley sum graphs and discover their perfect codes are discussed in Section \ref{Appen}.

\section{Subgroup Perfect Codes}\label{Results:SPC}
In this section, we begin by confirming a fundamental fact about perfect codes in Cayley sum graphs. We then present and prove a necessary and sufficient condition for a normal subgroup of a group to be, as well as not be, a subgroup perfect code of that group. Notably, we use these results to prove a statement that reveals the relationship between the subgroup perfect codes of a group and its quotient group. A similar statement for quotient groups was previously demonstrated in the context of Cayley graphs in \cite{ZZ}.
Prior to delving into the results, we recall
 a necessary and sufficient condition for a subgroup of a given group to be a subgroup perfect code.
\begin{lemma}\cite[Lemma~3.1]{Zhang-2024}\label{transversal}
Let $G$ be a group and $H$ a subgroup of $G$. Then $H$ is a subgroup perfect code of $G$ if and only if 
$G$ has a normal subset X such that $X \cup \{1\}$ is a left transversal of $H$ in $G$.
\end{lemma}
\begin{lemma}\label{lem:subgroup}
  Let $G$ be a group and $H$ a subgroup of $G$. Then $H$ is a subgroup perfect code of some Cayley sum graph of $G$ if and only if it is a subgroup perfect code of any subgroup of $G$ which contains $H$.
\begin{proof}
It is sufficient to prove the necessity. Suppose that $H \in \mathcal{P} (G, X)$, for some normal subset $X$ of $G$. By
 Lemma \ref{transversal}, $X\cup \{1\}$ is a left transversal
 of $H$ in $G$ which follows that every $g\in G\setminus H$ can be uniquely written as
 \begin{equation}\label{unique g}
 g = xh^{-1},
 \end{equation}
   for some $x\in X$ and $h\in H$. Let $K$ be  an arbitrary subgroup of $G$ such that $H\subseteq K$. From \eqref{unique g}, every $k\in K\setminus H$ can be uniquely written as $k=xh^{-1}$, for some $x\in X$ and $h\in H$. Since $x=kh$, we have $x\in X\cap K$. Therefore, 
    $(X\cap K)\cup \{1\}$ is a left transversal
 of $H$ in $K$. So, Lemma \ref{transversal} 
 implies that 
 $H \in \mathcal{P} (K, X\cap K)$.
\end{proof}
\end{lemma}
\begin{example}
Let $D_{2n}$ be the dihedral group of order $2n$ with  the following presentation
$$D_{2n}=\langle a,b:~ a^{n}=b^2=1, b^{-1}ab=a^{-1}\rangle.$$
Suppose that $n=4k+2$. In this case, the cyclic subgroup $\langle a^{2k+1}\rangle$ is 
a subgroup perfect code of $D_{2n}$,
see \cite[Example~5.8.]{Zhang-2024}. 
 So, every vertex in $D_{2n}$, including the vertices in $\langle a \rangle \setminus \langle a^{2k+1} \rangle$, is connected to either 1 or $a^{2k+1}$. Therefore, $\langle a^{2k+1} \rangle$ is also a perfect code of $\langle a \rangle$, as we expected from Lemma \ref{lem:subgroup}. 
 For $k=1$, you can see Figure \ref{fig:D_12} and Figure \ref{fig:SubD_12}.
\end{example}
\begin{figure}[ht]
 \centering
\begin{tikzpicture}[scale=1]
 \vertex[fill] (2) at (4,2) [label=above:\scriptsize{$a$}] {};
 \vertex[fill] (3) at (4,1) [label=above:\scriptsize{$a^5$}] {};
 \vertex[fill] (4) at (4,0) [label=above:\scriptsize{$ab$}] {};
 \vertex[fill] (5) at (4,-1) [label=above:\scriptsize{$a^3b$}] {};
 \vertex[fill] (6) at (4,-2) [label=above:\scriptsize{$a^5b$}] {};
 \vertex[fill] (7) at (8,2) [label=above:\scriptsize{$a^4$}] {};
 \vertex[fill] (8) at (8,1) [label=above:\scriptsize{$b$}] {};
 \vertex[fill] (9) at (8,0) [label=above:\scriptsize{$a^2b$}] {};
 \vertex[fill] (10) at (8,-1) [label=above:\scriptsize{$a^4b$}] {};
 \vertex[fill] (11) at (8,-2) [label=above:\scriptsize{$a^2$}] {};
 \tikzstyle{vertex}=[circle, draw, blue, inner sep=0pt, minimum size=4pt]
 \vertex[fill] (1) at (1,0) [ label=left:\scriptsize{$1$}] {};
 \vertex[fill] (12) at (11,0) [label=right:\scriptsize{$a^3$}] {};
 \path
 (1) edge (2)
 (1) edge (3)
 (1) edge (4)
 (1) edge (5)
 (1) edge (6)
 (12) edge (7)
 (12) edge (8)
 (12) edge (9)
 (12) edge (10)
 (12) edge (11)
 (2) edge (7)
 (2) edge (8)
 (2) edge (9)
 (2) edge (10)
 (3) edge (8)
 (3) edge (9)
 (3) edge (10)
 (3) edge (11)
 (4) edge (7)
 (4) edge (8)
 (4) edge (9)
 (4) edge (11)
 (5) edge (7)
 (5) edge (9)
 (5) edge (10)
 (5) edge (11)
 (6) edge (7)
 (6) edge (8)
 (6) edge (10)
 (6) edge (11);
 \end{tikzpicture}
\caption{$\mathrm{CS}(D_{12},\{a,a^5,ab,a^3b,a^5b\})$}
\label{fig:D_12}
\end{figure}
\begin{figure}[ht]
 \centering
\begin{tikzpicture}[scale=1]
 \vertex[fill] (2) at (3,1) [label=above:\scriptsize{$a$}] {};
 \vertex[fill] (3) at (3,-1) [label=above:\scriptsize{$a^5$}] {};
 \vertex[fill] (4) at (5,1) [label=above:\scriptsize{$a^2$}] {};
 \vertex[fill] (5) at (5,-1) [label=above:\scriptsize{$a^4$}] {};
 \tikzstyle{vertex}=[circle, draw,blue, inner sep=0pt, minimum size=4pt]
 \vertex[fill] (1) at (2,0) [ label=left:\scriptsize{$1$}] {};
 \vertex[fill] (6) at (6,0) [label=right:\scriptsize{$a^3$}] {};
 \path
 (1) edge (2)
 (1) edge (3)
 (2) edge (5)
 (3) edge (4)
 (6) edge (4)
 (6) edge (5);
 \end{tikzpicture}
\caption{$\mathrm{CS}(\langle a \rangle ,\{a,a^5\})$}
\label{fig:SubD_12}
\end{figure}
We recall that the \textit{center} of a group $G$, denoted $Z(G)$, is defined as the set of elements in $G$ that commute with every element in $G$. Additionally, for all $x, g \in G$, the  \textit{commutator} of $x$ and $g$ is defined as $[x,g] = x^{-1}g^{-1}xg$. 
Furthermore, for $x, g \in G$, the conjugate of $x$ by $g$ is given by $x^g := g^{-1}xg$. For $S \subseteq G$, the conjugate of $S$ by $g$ is $S^g = \{s^g : s \in S\}$.
\begin{theorem}\label{xh}
Let $G$ be a group and $H$ a normal subgroup of $G$. Then 
$H$ is a subgroup perfect code of some Cayley sum graph of $G$ if and only if
for all $x\in G$, $[x, G] \subseteq H$ implies  that $xh\in Z(G)$ for some $h \in H$.
\begin{proof}
Suppose that $H\in \mathcal{P}(G,Y)$, where $Y$ is a normal subset of $G$. Then from Lemma \ref{transversal}, $L=Y\cup \{1\}$ is a left transversal of $H$ in $G$. Suppose that $ x\in G$ and $[x, G] \subseteq H$.
Since $L$ is a left transversal of $H$ in $G$, there exists an element $y\in L$ such that $x\in yH$ and so, $xH=yH$. Then we have
\begin{equation}\label{xH=yH}
x^gH=Hx^g=Hy^g=y^gH, ~~ \forall g\in G,
\end{equation}
because $H$ is normal in $G$. Moreover, $[x, G] \subseteq H$ implies that
$x^{-1}x^g  H=H$ and $x^g  H=xH$, for all $g\in G$. From \eqref{xH=yH}, it is concluded that $y^g  H=yH$ and consequently,  $y=y^g$ since $y\in L$ and $y^g \in  L^g = L$. So, from the fact that $x\in yH$,  there exists an element
$h\in H$ for which $xh=(xh)^g$.

Suppose that for any $x\in G$, $[x, G] \subseteq H$ implies $xh\in Z(G)$ for some $h\in H$.
Consider the set $I:=\{y_1H,\ldots,y_sH\}\subseteq G/H$ such that $y_iH=
y_i^gH$ where $g\in G$ and $1\leqslant i\leqslant s$. Also, set $J:=\{g_1H, g_1^gH,\ldots,g_tH, g_t^gH\}\subseteq G/H$ such that $g_iH\neq
g_i^gH$ where $g\in G$ and $1\leqslant i\leqslant t$. For any element $y_iH \in I$,  we have $  y_i^{-1}y_i^gH\subseteq H$ and so, $ y_i^{-1}y_i^g=[y_i,g]\in H$. Therefore, there exists $h_i\in H$ such that $y_ih_i \in Z(G)$ and $y_ih_i=(y_ih_i)^g$.  Denote $z_i = y_ih_i$ and set $L=\{z_1,\ldots,z_s,g_1,g_1^g,\ldots,g_t,g_t^g\}$.
We can see that $L\cup \{1\}$ is a left transversal of $H$ in $G$ and $L^g=L$. From Lemma \ref{transversal}, we have $H\in \mathcal{P}(G,L)$.  
\end{proof}
\end{theorem}
\begin{remark}
Theorem \ref{xh} implies that every subgroup of  an abelian group is a perfect code of some Cayley sum graph, as proved in \cite[Theorem~5.1]{Zhang-2024}.
\end{remark}
\begin{theorem}\label{xH}
Let $G$ be a group and $H$ a normal subgroup of $G$. Then $H$ is not a subgroup perfect code of any Cayley sum graphs of $G$ if and only if there exists a left coset $L = xH$ with $L = L^{g}$ (that is, $g^{-1}Lg=L$ for all $g\in G$)  which contains no element of $Z(G)$.
In particular, if $H$ is not a subgroup perfect code of any Cayley sum graphs of $G$, then there exists an element $x \in G \setminus H$ such that $[x, G]  \subseteq  H$, and $xH$ contains no element of $Z(G)$.

\begin{proof}
Suppose that there is a left coset $L = xH$ with $L = L ^{g}$, where $g\in G$,
 containing no element of $Z(G)$. 
 Assume the contrary that $H\in \mathcal{P} (G, X)$, where $X$ is a normal subset of $G$. By Lemma \ref{transversal}, $T=X\cup \{1\}$ is a left transversal of $H$ in $G$. Hence, 
 there exists a unique element $y\in T$ such that $T \cap  L = \{y\}$.  Since both $T$ and $L$ are
normal subsets of $G$,  we have 
$$
\{y\}^g= (T\cap L)^g
=T^g\cap L^g=T\cap L=\{y\},$$
for all $g\in G$.  Therefore, we have $y^g=y$, for all $g\in G$ which implies that $y \in Z(G)$.  Hence, $L$ includes at least one element of $Z(G)$ which contradicts our assumption.

Now assume that $H$ is not a perfect code of $G$.
Take a subset $T$ of $G$ with maximum cardinality such that 
\begin{itemize}
\item [(i)] $1 \in T,$
\item [(ii)] $T^g = T,~~ \forall g\in G,$
\item [(iv)] $xH\neq yH, ~~\forall x, y\in T ~~ s.t. ~~ x\neq y.$
\end{itemize}
The existence of $T$ follows from the fact that there are subsets of $G$, say $\{1\}$, with all these properties. Since $H$ is not a perfect code of $G$, by Lemma \ref{transversal}, $T$ is not a left transversal of $H$ in $G$. It follows that $G \setminus TH \neq \emptyset$ and therefore, we can take an element $x\in G \setminus TH$. Set $L = xH$. Then, we have $L = x^gH$, where $g\in G$, because otherwise $T\cup \{x, x^g\}$
satisfies all properties $(i)-(iii)$ and so, contradicts the maximality of $T$. 
Also, since $H\trianglelefteq G$, we have
$$L^{g}=(xH)^g=g^{-1}xHg=g^{-1}xgH=x^{g}H=L.$$

Now we prove that $L$ contains no element of $Z(G)$. Suppose to the contrary and let $z\in L\cap Z(G)$. Then set $X=T\cup \{z\}$. It is straightforward to check that  all properties $(i)-(iii)$ hold for $X$ and so,
contradicts the maximality of $T$. 


Now we prove the last statement in the theorem. Assume that $H$ is not a perfect code of any  Cayley sum graph of $G$. Then by what we have proved above there exists $y\in G\setminus H$ such that $yH = y^{g}H$,
where $g\in G$. So, we have $y^{g}= yh$ for some $h \in H$. Set 
$x=yh$. We have
$$x^{-1}x^g=h^{-1}y^{-1}y^gh^g=h^{-1}y^{-1}yhh^g=h^g,$$
and since $H\trianglelefteq G$, $x^{-1}x^g\in H$. Furthermore,
$xH=yH$ and so, contains no element of $Z(G)$.
\end{proof}
\end{theorem}
\begin{remark}
It is worth noting that by considering the negation of the statement in Theorem \ref{xh}, we can obtain the last statement provided in Theorem \ref{xH}.
\end{remark}
\begin{theorem}\label{Qutient}
Let $G$ be a group and $N$ and $H$  normal subgroups of $G$ such that $H$  contains $N$. Then the following hold:
\begin{itemize}
\item[(i)] if $H$ is a subgroup perfect code of some Cayley sum graph of $G$, then $H/N$ is a subgroup perfect code of some Cayley sum graph of $G/N$;
\item[(ii)] if $N$ and $H/N$ are subgroup perfect codes of some Cayley sum graphs of $G$ and $G/N$, respectively, then $H$ is a subgroup perfect code of some Cayley sum graph of $G$.
\end{itemize}
\begin{proof}
(i). Suppose that $H \in \mathcal{P} (G, X)$, where $X$ is a normal subset of $G$. Therefore,  we have
\begin{equation}\label{g-G-H}
\forall g\in G\backslash H,~~ \exists ! ~ h\in H ~~~s.t. ~~~gh=hg\in X.
\end{equation}
Now we define the set 
$X/N:=\{xN:~ x\in X\}$ with the same binary and inverse operations  as in quotient groups. It is straightforward to see that $X/N$ is a normal subset of $G/N$. For every $gN\in G/N\setminus H/N$, it is concluded from \eqref{g-G-H} that there exists a unique element 
$hN\in H/N$ such that $ghN\in X/N$ and hence,  $H/N \in \mathcal{P} (G/N, X/N)$.

(ii). Suppose to the contrary that $H$ is not a perfect code of any Cayley sum graphs of $G$. By Theorem \ref{xH}, there exists a element $x \in G\setminus H$ such that $[x, G]\subseteq H$ and $xH$ contains no element of $Z(G)$.  Since $N$ is a normal subgroup of $G$ contained in $H$,
$xN\in G/N\setminus H/N$ and $x^{-1}x^gN\in H/N$, for all $g\in G$. So,
we have $x^gN=xN\in H/N$ and subsequently, $xN (H/N)=x^gN(H/N)$.
Since $H/N$ is a perfect code of $G/N$, by Theorem \ref{xH}, 
$xN(H/N)$ contains at least one element of $Z(G)$. So, there exists a central element $aN \in H/N$ such that $(xNaN)^g= xNaN$ and so, $(xa)^{-1}(xa)^gN=N$. Hence, $(xa)^{-1}(xa)^g\in N$.
 Since $N$ is a perfect code and a normal subgroup of $G$, by Theorem \ref{g-G-H}, there exists $b \in N$ such that $(xab)^g=xab$. Note that $a, b\in H$ and $xab\in xH \cap Z(G)$. However, $xH$ contains no element of $Z(G)$, a contradiction.
\end{proof}
\end{theorem}
  
\section{Subgroup Perfect Codes and Total Perfect Codes of  $V_{8n}$}\label{special}
In \cite{Zhang-2024}, perfect codes have been explored within the context of Cayley sum graphs for certain groups, such as abelian and dihedral groups. 
Here, the focus is on exploring the perfect codes of the group $V_{8n}$, which is one of the famous and applicable groups. Before proceeding, we recall a lemma pertaining to the connectivity of a Cayley sum graph in \cite{AT} and a necessary and sufficient condition for total perfect codes proven in \cite{Zhang-2024}.  Note that, although the definition of Cayley sum graphs in \cite{AT} differs slightly from ours, the aforementioned result still holds for our graphs. This is because our Cayley sum graphs are obtained simply by removing loops from the graphs in \cite{AT}, preserving the essential properties relevant to the stated result. 
\begin{lemma}\cite[Theorem~2]{AT}\label{connect}
 Let $X$ be a normal subset of a group $G$. Then the Cayley sum graph $\rm{CS}(G, X)$ is connected if and only if $G = \langle X \rangle$ and $|G : X^{-1}X|\leqslant 2$.
\end{lemma}

\begin{lemma}\label{Total-Zhang-2024}\cite[Lemma 4.1.]{Zhang-2024}
Let $G$ be a group and $H$ a subgroup of $G$. Then H is a subgroup total perfect code of some Cayley sum graph of $G$ if and only if there exists a normal subset $Y$ of $G$  such that $Y$ is a left transversal of $H$ in $G$ and $H\cap Y$ contains precisely one unique element which is a nonsquare of $H$.
\end{lemma}
We recall that for an odd positive integer $n$, the group $V_{8n}$  defined in  \cite[p. 178]{James1993} is as follows 
$$V_{8n}=\langle a,b:~ a^{2n}=b^4=1, ba=a^{-1}b^{-1}, b^{-1}a=a^{-1}b\rangle.$$ 
\begin{example}\label{Gamma0}
 We introduce the set 
 $$X:=\{b, ab, a^2b, a^3b, \dots, a^{2n-1}b,b^3, ab^3,  a^3b^3, a^3b^3, \dots, a^{2n-1}b^3\},$$
  which is normal in $V_{8n}$ according to the conjugacy classes detailed in Section \ref{Appen}. The Cayley sum graph $\Gamma_0:=\mathrm{CS}(X,V_{8n})$ is isomorphic to 
 the complete bipartite graph $K_{4n,4n}$ with two parts $X$ and $$S:=V_{8n} \setminus X=\{e, a, a^2, a^3, \dots, a^{2n-1},b^2, ab^2,  a^3b^2, a^3b^2, \dots, a^{2n-1}b^2\}.$$ 
 Also,  the subgroups $\langle a^jb\rangle$ and $\langle a^jb^3\rangle$, where $j$ is an odd number in the range $1 \leqslant j < 2n$, are total perfect codes of $\Gamma_0$.
 
\end{example}
\begin{example}\label{Gamma1-Gamma2} 
 By considering the conjugacy classes of $V_{8n}$ discussed in Section \ref{Appen}, the following sets are normal in $V_{8n}$: 
  \begin{align*} 
  &Y_1:=\{ab, a^3b, \dots, a^{2n-1}b, ab^3, a^3b^3, \dots, a^{2n-1}b^3\},\\
 & Y_2:=\{b, a^2b, a^4b, \dots, a^{2n-2}b, b^3, a^2b^3, a^4b^3, \dots, a^{2n-2}b^3\},\\
&Z:=\{a, a^3, \dots, a^{2n-1}, ab^2, a^3b^2, \dots, a^{2n-1}b^2\},\\
&Z^{\prime}:=\{a^2, a^4, \dots, a^{2n-2}, a^2b^2, a^4b^2, \dots, a^{2n-2}b^2\}\cup \{b^2\}.
\end{align*}
Using these sets, we define four Cayley sum graphs as $\Gamma_1:=\mathrm{CS}(V_{8n},Y_1\cup Z)$,  $\Gamma_1^{\prime}:=\mathrm{CS}(V_{8n},Y_1\cup Z^{\prime})$,  $\Gamma_2 := \mathrm{CS}(V_{8n}, Y_2\cup Z)$ and  $\Gamma_2^{\prime} := \mathrm{CS}(V_{8n}, Y_2\cup Z^{\prime})$. It can be easily verified that the subgroups $\langle a^n\rangle$, $\langle a^nb^2\rangle$, $\langle a^jb\rangle$ and $\langle a^jb^3\rangle$, where $1\leqslant j<2n$ is odd, are total perfect codes of $\Gamma_1$. Furthermore, $\langle a^n\rangle$ and $\langle a^nb^2\rangle$ are perfect codes, while $\langle a^jb\rangle$ and $\langle a^jb^3\rangle$, where $1\leqslant j<2n$ is odd, are total perfect codes of $\Gamma_1^{\prime}$. It is also straightforward to check that the subgroups $\langle a^n\rangle$ and $\langle a^nb^2\rangle$ are total perfect codes, while $\langle a^jb\rangle$ and $\langle a^jb^3\rangle$, where $1\leqslant j<2n$ is odd, are perfect codes of $\Gamma_2$. Moreover, $\langle a^n\rangle$, $\langle a^nb^2\rangle$, $\langle a^jb\rangle$ and $\langle a^jb^3\rangle$, where $1\leqslant j<2n$ is odd, are  perfect codes of $\Gamma_2^{\prime}$. 
 \end{example}
We mention that in the following theorems, we consider $j$ to be an odd number in the range $1 \leqslant j < 2n$ and use the normal sets $Z$, $Z^{\prime}$, $Y_1$, and $Y_2$ described in Example \ref{Gamma1-Gamma2}. 
\begin{theorem}\label{Th:V_8nPerfect}
For an odd positive integer $n$, the nontrivial subgroup perfect codes in any connected Cayley sum graph of $V_{8n}$ are limited to the subgroups $\langle a^n\rangle$, $\langle a^nb^2\rangle$, $\langle a^jb\rangle$, and $\langle a^jb^3\rangle$. In other words, there exist exactly $2n+2$ subgroups of $V_{8n}$ that can function as nontrivial perfect codes in some connected Cayley sum graph of $V_{8n}$.
\begin{proof}
Suppose that $\Gamma: =\mathrm{CS} (V_{8n},X)$ is a connected Cayley sum graph. We know that $X$ is a normal subset of $V_{8n}$ and  Lemma \ref{connect} implies that $\langle X \rangle = V_{8n}$.
So, we have $X\cap \langle a\rangle b\neq \emptyset$. Let $H$ be 
a nontrivial subgroup perfect  code of $\Gamma$. From Lemma \ref{transversal}, $X\cup \{1\}$ is a left transversal of $H$ in 
$V_{8n}$ and therefore, $|H|\left(|X|+1\right)=8n$. Since $|H|>1$, we conclude that $|X|<4n$. Hence, $X$ contains exactly one of the following two conjugacy classes $Y_1$ and $Y_2$.
  Set $Y:=X\cap \langle a\rangle b$. We have $|Y|=2n$ which leads to
$2n<|X|<4n$. Consequently, $|H|=2$ or 3. We obtain from \cite[Section~2.4]{SA2021} that  $H$ is of one of the following forms:
$$
\begin{array}{ll}
(i)~ H=\langle a^{\frac{2n}{3}}\rangle, \text{where}~ 3|2n,
&(iv)~ H=\langle b^2\rangle\\
(ii)~ H=\langle a^n\rangle, &(v)~ H=\langle a^jb\rangle, \\
  (iii)~ H=\langle a^nb^2\rangle, & (vi)~ H=\langle a^jb^3\rangle,  
\end{array}
$$
Noticeably, the case $(i)$ does not happen because it is equivalent to  $|X|+1=8n/3$. This implies that $|X|<|Y|=2n$  which 
is a contradiction.
To continue  our discussion, we consider the following two cases that pertain to the set $Y$.\\
$\mathbf{Case ~1:}$ ~ $Y=Y_1$. We analyze the various outcomes when considering different forms of $H$.
\begin{itemize}
\item  $H=\langle a^n\rangle$ or $H=\langle a^n b^2\rangle$.  When considering the left cosets of $H$ in \eqref{Cosets1} and \eqref{Cosets2},  we observe that the sets  $Z$ and $Z^{\prime}$ are the unique sets for which
$Y_1\cup Z$ and $Y_1\cup Z^{\prime} \cup \{1\}$  normal in $V_{8n}$ and form left transversals of
$H$ in $V_{8n}$. 
\item $H=\langle b^2\rangle$.  By employing \eqref{Cosets3} in this case, it becomes impossible to construct a left transversal of $H$ in $G$ using a set with cardinality less than or equal to $2n$.
\item $H=\langle a^j b\rangle$ or $H=\langle a^j b^3\rangle$. 
Upon analyzing the left cosets of $H$ in \eqref{Cosets4} and \eqref{Cosets5}, we observe that the set $Z$ is the unique set  for which  $Y_1\cup Z$ is normal and a left transversal of $H$ in $V_{8n}$.
\end{itemize}
$\mathbf{Case~ 2:} ~Y=Y_2$.
 Now we discuss different forms of  $H$. 
\begin{itemize}
\item $H=\langle a^n\rangle$ or $H=\langle a^nb^2\rangle$. From the left cosets described in (\ref{Cosets1}) and \eqref{Cosets2}, 
$Y_2\cup Z$ and $Y_2\cup Z^{\prime}\cup \{1\}$ are the only possible normal subsets of $V_{8n}$ that can be considered as left transversals of $H$ in $V_{8n}$.
\item $H=\langle b^2\rangle$. By considering \eqref{Cosets3}, it is impossible to identify a normal subset of $V_{8n}$ with cardinality less than or equal to $2n$ whose union with $Y_2$ forms a left transversal of $H$ in $V_{8n}$.
\item $H=\langle a^j b\rangle$ or  $H=\langle a^j b^3\rangle$. By \eqref{Cosets4} and \eqref{Cosets5}, 
$Y_2\cup Z^{\prime}\cup \{1\}$ is the only left transversal of $H$ in $V_{8n}$ that can be considered.
\end{itemize}
Therefore, $\Gamma$ can be one of the graphs described in Example \ref{Gamma1-Gamma2}.
Note that the two cases, Case 1 and Case 2, encompass all possible scenarios for connected Cayley sum graphs of $V_{8n}$. As a result, they cover all subgroup perfect codes. This is due to the requirement that, regardless of the chosen conjugacy classes to be contained in $X$, $X$ must include one of the conjugacy classes from either Case 1 or Case 2. This is because $X$ needs to generate $V_{8n}$, and none of the other conjugacy classes can generate the elements of these two sets.
\end{proof}
\end{theorem}
\begin{theorem}\label{Th:V_8n}
For an odd positive integer $n$, the nontrivial subgroup total perfect codes in any connected Cayley sum graph of $V_{8n}$ are limited to the subgroups $\langle a^n\rangle$, $\langle a^nb^2\rangle$, $\langle a^jb\rangle$, and $\langle a^jb^3\rangle$. In other words, there exist exactly $2n+2$ subgroups of $V_{8n}$ that can function as nontrivial total perfect codes in some connected Cayley sum graph of $V_{8n}$.
\begin{proof}
Suppose that $\Gamma: =\mathrm{CS} (V_{8n},X)$ is a connected Cayley sum graph. We know that $X$ is a normal subset of $V_{8n}$ and  Lemma \ref{connect} implies that $\langle X \rangle = V_{8n}$.
So, we have $X\cap \langle a\rangle b\neq \emptyset$. Let $H$ be 
a nontrivial subgroup total perfect  code of $\Gamma$. From Lemma \ref{Total-Zhang-2024}, $X$ is a left transversal of $H$ in 
$V_{8n}$ and therefore, $|H||X|=8n$. Since $|H|>1$, we have $|X|\leqslant 4n$.
Consequently, $X$ must contain either the union $Y_1 \cup Y_2$ or exactly one of the two conjugacy classes $Y_1$ and $Y_2$. In the case where $Y_1 \cup Y_2 \subseteq X$, it follows that $X = Y_1 \cup Y_2$ due to $|Y_1 \cup Y_2| = 4n$. As a result, $\Gamma$ is equivalent to $\Gamma_0$ as presented in Example \ref{Gamma0}. If either $Y_1$ or $Y_2$ is a subset of $X$, the remaining proof closely mirrors the argument detailed in Theorem \ref{Th:V_8nPerfect}.

\end{proof}
\end{theorem}
\section{Appendix: Computational Details}\label{Appen}
 In the appendix, readers will find computational details that explore the relationships between the members of the group $V_{8n}$, the conjugacy classes of $V_{8n}$, and the left cosets of a subgroup $H$ within $V_{8n}$. These details played a crucial role in the construction of the perfect code in the group $V_{8n}$.\\
 \textbf{Relationships between elements of $V_{8n}$:}

 Within the group of $V_{8n}$, the following relationships hold
$$
a^ib^j=
\begin{cases}
b^{-j}a^i~~~~~~\text{if $i$ is odd and $j$ is even },\\
b^{-j}a^{-i}~~~~\text{if $i$ and $j$ are odd},\\
b^{j}a^i~~~~~~~~\text{if $i$ and $j$ are even },\\
b^{j}a^{-i}~~~~~~\text{if $i$ is even and $j$ is odd}.\\
\end{cases}
$$
 \textbf{Conjugacy Classes:}
 
Identifying normal sets within the group $V_{8n}$ required knowledge of the group's conjugacy classes. The conjugacy classes for $V_{8n}$ are as follows
\begin{align*}
&\{1\}, \{b^2\}, \{a^{2r+1}, a^{-2r-1}b^2\},~ r=0,\dots, n-1,\\
& \{a^{2s}, a^{-2s}\}, \{a^{2s}b^2, a^{-2s}b^2\},~ s=1,\dots, \frac{n-1}{2},\\
&\{a^jb^k: j ~\text{even},~ k=1~ \text{or}~ 3\},\\
&\{a^jb^k: j~ \text{odd}, ~k=1 ~\text{or} ~3\},
\end{align*}
see \cite[p. 421]{James1993}.\\
 \textbf{Left Cosets and Left Transversal:}
 
  To determine a left transversal of a subgroup $H$ in $G$, we calculated the left cosets of $H$ in $G$ for five different cases of $H$:
  \begin{itemize}
  \item $H=\langle a^n\rangle$. In this case, the left cosets of $H$ are
\begin{align}\label{Cosets1}
&a^iH=\{a^i, a^{n+i}\}, a^ibH=\{a^ib, a^{i-n}b^3\},  a^ib^2H=\{a^ib^2, a^{i+n}b^2\}, a^ib^3H=\{a^{i-n}b, a^ib^3\},\nonumber\\
&bH=\{b, a^{-n}b^3\}, b^2H=\{b^2, a^{n}b^2\}, 
 b^3H=\{b^3, a^{-n}b\},
\end{align}
where $1\leqslant i \leqslant 2n-1$.
\item $H=\langle a^n b^2\rangle$. The left cosets of $H$ are
\begin{align}\label{Cosets2}
&a^iH=\{a^i, a^{n+i}b^2\}, a^ibH=\{a^ib, a^{i-n}b\},  a^ib^2H=\{a^{i+n},a^ib^2\}, a^ib^3H=\{a^{i}b^3, a^{i-n}b^3\},\nonumber\\
&bH=\{b, a^{-n}b\}, b^2H=\{b^2, a^{n}\}, 
 b^3H=\{b^3, a^{-n}b^3\},
\end{align}
where $1\leqslant i \leqslant 2n-1$. 
\item $H=\langle b^2\rangle$. The left cosets of $H$ are
\begin{align}\label{Cosets3}
&a^iH=\{a^i, a^{i}b^2\}, a^ibH=a^ib^3H=\{a^ib, a^{i}b^3\}, \nonumber\\
& a^ib^2H=\{a^{i},a^ib^2\}, bH=b^3H=\{b, b^3\},
\end{align}
where $1\leqslant i \leqslant 2n-1$. 
\item $H=\langle a^j b\rangle$, where $1\leqslant j<2n$ is odd. In this case, the left cosets of $H$ are
\begin{align}\label{Cosets4}
&a^iH=\{a^i, a^{i+j}b\}, a^ibH=\{a^{i-j}, a^{i}b\},  a^ib^2H=\{a^ib^2, a^{i+j}b^3\}, a^ib^3H=\{a^{i-j}b^2, a^ib^3\},\nonumber\\
&bH=\{a^{-j}, b\}, b^2H=\{b^2, a^{j}b^3\}, 
 b^3H=\{b^3, a^{-j}b^2\},
\end{align}
where $1\leqslant i \leqslant 2n-1$. 
\item $H=\langle a^j b^3\rangle$, where $1\leqslant j<2n$ is odd. In this case, the left cosets of $H$ are
\begin{align}\label{Cosets5}
&a^iH=\{a^i, a^{i+j}b^3\}, a^ibH=\{a^{i}b, a^{i-j}b^2\},  a^ib^2H=\{a^{i+j}b, a^{i}b^2\}, a^ib^3H=\{a^{i-j}, a^ib^3\},
\nonumber\\
&bH=\{b, a^{-j} b^2\}, b^2H=\{b^2, a^{j}b\}, 
 b^3H=\{ a^{-j}, b^3\},
\end{align}
where $1\leqslant i \leqslant 2n-1$.
  \end{itemize}
  \section*{Acknowledgement}
  The first author was supported by University of Mazandaran, Grant Number 60673.

\end{document}